\documentclass[preprint]{elsarticle}
\usepackage[utf8]{inputenc}
\usepackage{bm}
\usepackage{amsmath}
\usepackage{amsfonts}
\usepackage{amssymb}
\usepackage{graphicx}
\usepackage{xcolor}
\usepackage{cancel}
\usepackage{caption}
\usepackage{subfigure}
\usepackage[ruled,vlined]{algorithm2e}

\usepackage{tikz}
\usepackage{pgfplots}

\newcommand*{\Scale}[2][4]{\scalebox{#1}{$#2$}}%
\newcommand{\half}{\Scale[0.6]{\tfrac{1}{2}}}

\usepackage[T1]{fontenc}
\usepackage{lmodern}
\usepackage[top=3cm,bottom=4cm,left=3cm,right=3.2cm,asymmetric]{geometry}

\newcommand{\N}[1]{\check{#1}}
\newcommand{\D}[1]{\bar{#1}}

\SetKwProg{Fn}{Function}{}{}
\DontPrintSemicolon
\SetKwInOut{Output}{output}
\SetKwInOut{Input}{input}
\SetKw{KwStep}{ step }

\begin{document}
\begin{frontmatter}

\title{A spectral-Galerkin turbulent channel flow solver for large-scale simulations}
\author[mmo]{Mikael Mortensen}
\ead{mikaem@math.uio.no}
\address[mmo]{Department of Mathematics, Division of Mechanics, University of Oslo}

\begin{abstract}
A fully (pseudo-)spectral solver for direct numerical simulations of large-scale turbulent channel flows is described. The solver utilizes the Chebyshev base functions suggested by J. Shen [SIAM J. Sci. Comput., 16, 1, 1995], that lead to stable and robust numerical schemes, even at very large scale. New and fast algorithms for the direct solution of the linear systems are devised, and algorithms and matrices for all required scalar products and transforms are provided. We validate the solver for very high Reynolds numbers. Specifically, the solver is shown to reproduce the first order statistics of Hoyas and Jim\'{e}nez [Phys. Fluids, 18(1), 2006], for a channel flow at $Re_{\tau}=2000$. The solver is available through the open source project spectralDNS [https://github.com/spectralDNS].
\end{abstract}

\end{frontmatter}
\section{Introduction}
Direct Numerical Simulations (DNS) of turbulent flows is a very important research tool, utilized across a range of scientific communities \cite{Moin98}. DNS is used extensively to validate statistical models, and to further our understanding of complex mechanisms taking place inside turbulent flows. One of the many advantages of DNS is that it provides all information about a flow, and quantities that can be very hard to study experimentally, like velocity-pressure interactions, are trivially extracted from a DNS. In this regard, DNS both complements and extends the knowledge we are able to extract from experiments.

The most commonly known DNS use simple geometries, because turbulence physics may then most easily be isolated and studied. Isotropic and homogeneous turbulence are usually studied numerically in triply periodic domains, which allows for a spectral Fourier decomposition in all three spatial directions. Spectral methods are often favored in DNS due to their superior accuracy and resolution properties. One example is given in the DNS review of Moin and Mahesh \cite{Moin98}, who report that, for similar accuracy in first derivatives, a second-order finite difference scheme requires approximately 5.5 times more points than Fourier, in each spatial direction, whereas for a 6'th order Pad\'{e} scheme the factor is about 1.6.

In this paper we will consider the pressure driven turbulent channel flow, where there are two periodic directions that can be handled with Fourier expansions, and a non-periodic wall-normal direction that requires a different type of discretization. There are many challenges associated with this inhomogeneity not faced by the pure Fourier solvers, but the first problem at hand is the discretization. Early DNS channel solvers, see, e.g., \cite{Moin80, Kleiser80, Kim87}, typically used a Chebyshev expansion for the wall-normal direction and, as such, were still able to obtain spectral accuracy in all three spatial directions. A Chebyshev-tau technique, that utilize the recurrence relations of the Chebyshev polynomials, was used to approximate derivatives, and the coefficient matrices that appeared (tridiagonal) could then be inverted directly and efficiently \cite{Kim87}. A downside to the Chebyshev-tau method is usually quoted \cite{canuto1988} as numerical instability and roundoff errors, caused by the recurrence relation, and severe condition numbers of the coefficient matrices. For Chebyshev-tau methods the condition numbers have been reported to grow with size as $\mathcal{O}(N^8)$, for a discretization using $N$ points in the wall-normal direction. Discouraged by these numbers, all major recent channel flow simulations have found other, non-spectral, ways of discretizing the non-periodic direction.

The largest known channel simulations to date have been performed by Lee and Moser \cite{leemoser15}, where, for $Re_{\tau}\approx 5200$, they used a computational box of resolution $[10240 \times 1530 \times 7680]$ for the streamwise, wall-normal and spanwise directions, respectively. Lee and Moser used seventh-order B-splines for the wall-normal direction. Other simulations of similar magnitude have been performed by Hoyas and Jim\'{e}nez \cite{hoyas06, hoyas08}, Lozano-Duran and Jim\'{e}nez \cite{Lozano2014} and Bernardini, Pirozzoli and Orlandi \cite{bernardini2014}. Bernardini et al. used second-order finite differences throughout. The Jim\'{e}nez group used dealiased Fourier in the two periodic directions and seventh-order compact finite differences, with fourth-order consistency and extended spectral-like resolution \cite{Lele92}, for the wall-normal direction. The solver by Jim\'{e}nez' group is reported to switch from Chebyshev to finite differences if the resolution is above a certain threshold \cite{hoyas08} (reached around $Re_{\tau}=1000$). In other words, they attempt to use a fully spectral discretization for as large $Re_{\tau}$ as possible. As previously mentioned, spectral methods are attractive for their accuracy and resolution properties, that are superior to those of any finite difference or spline method. As such, it is desirable to develop fully spectral solvers that can be used for large-scale turbulence simulations.

In his seminal papers \cite{Shen94,Shen95}, Jie Shen describe how to construct Legendre and Chebyshev basis functions that lead to sparse matrices, susceptible to very fast direct solvers. To the author's knowledge, the bases have not been used for large-scale channel flow simulations, and algorithms for the required direct solvers have not, until now, been devised. Shen's bases have been used for the Navier Stokes equations before, though. Bouchon et.al.  \cite{Bouchon01} describe a spectral-Galerkin formulation very similar to the one used in this paper. However, they choose the Legendre basis over Chebyshev, which has some consequences when aiming for large-scale, because fast transforms are required in moving from spectral to physical space, and back again. For Fourier and Chebyshev bases, the Fast Fourier Transforms (FFTs of $\mathcal{O} (N \log N)$) apply directly. However, until recently, the discrete Legendre transforms required $\mathcal{O}(N^2)$ operations. This scaling has now been improved by several authors, as recently reviewed by Hale and Townsend \cite{Hale14, Hale2015}, but the methods are still not quite on par with the FFTs. For example, Hale and Townsend describe an $\mathcal{O}(N (\log N)^2 / \log \log N)$ algorithm, using intermediate fast transforms from Legendre to Chebyshev coefficients.

In this paper we will describe and assess a spectral-Galerkin channel flow solver based on Fourier and Shen's Chebyshev basis \cite{Shen95}. The solver will consist of parts that scale at worst as $N \log N$, for a 1D problem of size $N$, and as such as $N^3 \log N$ for a 3D box of size $N^3$. We will give a proper description of the theoretical basis in Sec~\ref{sec:prelim}, the discretization of Navier-Stokes equations in Sec~\ref{sec:discretizationNS}, and we will describe necessary algorithms, including a new fast direct solver for the biharmonic problem that arise, in Sec~\ref{sec:implementation}. We will finally show, in Sec~\ref{sec:verification}, that roundoff is not a major issue, and that the Shen-Fourier spectral-Galerkin method is in deed applicable to large-scale simulations. The solver has been implemented in the open source code spectralDNS \cite{spectralDNS}, where the bulk of the code is written in high-level Python \cite{python}, with critical parts migrated to Cython \cite{cython} for efficiency.

\section{Basis functions and fast transforms}
\label{sec:prelim}
The Navier-Stokes equations, used to describe turbulent flow in a doubly 
periodic channel, can be written in rotational form as
\begin{align}
 \frac{\partial \bm{u}}{\partial t}   &= \bm{\mathcal{H}} + \nu 
 \nabla^2 \bm{u} - \nabla{\tilde{p}}, \notag \\
 \nabla \cdot \bm{u} &= 0, \label{eq:NS}
\end{align}
where $\bm{u}(\bm{x}, t)=(u, v, w)$ is the velocity vector, $\bm{x}=(x, y, z)$ 
and $t$ are position and time, and the nonlinearity $\bm{\mathcal{H}}(\bm{x}, t) = (\mathcal{H}_x, \mathcal{H}_y, \mathcal{H}_z) = \bm{u}\times 
\bm{\omega}$, where $\bm{\omega} = \nabla \times \bm{u}$.  The constant dynamic viscosity is denoted as $\nu$ 
and $\tilde{p}(\bm{x}, t)$ is a pressure modified to account for both the driving force, $\beta(t)$, and the kinetic energy, i.e., $\tilde{p} = p + \bm{u} \cdot \bm{u}/2 + \beta y$, where $p$ is the instantaneous 
pressure normalized by a constant density. The computational domain is 
$\Omega=[-1, 1]\times [0, L_y] \times [0, L_z]$, with channel walls 
located at $x=\pm 1$, such that no-slip applies at 
$ \bm{u}(\pm 1, y, z, t) = 0$. The walls are spanning 
the $y-z$ plane and the equations are periodic in the $y$ and $z$ directions with periodic lengths $L_y$ and $L_z$, respectively.

The domain $\Omega$ is discretized using $N = (N_x, N_y, N_z)$ intervals, where the two periodic directions use uniform intervals. The computational mesh is given as
\begin{align}
X_N = \Big\{ \bm{x} \in \mathbb{R}^3 | &(x_i, y_j, z_k) = \left( h(i), \frac{jL_y}{N_y}, \frac{kL_z}{N_z} \right), \text{where} \notag \\
&(i, j, k) \in [0, 1, \ldots, N_x] \times [0, 1, \ldots, N_y-1] \times [0, 1, \ldots, N_z-1] \Big\}, \label{eq:Xn}
\end{align}
where $x_i = h(i)$ represents
\begin{equation}
h(i) = \begin{cases}
\cos \left(\frac{i \pi }{N_x} \right) \, &\forall \, i=0,1, \ldots, N_x \quad  \text{for Chebyshev-Gauss-Lobatto points}, \\
\cos \left(\frac{(2i +1)\pi}{2N_x+2} \right) \, &\forall \, i=0,1, \ldots, N_x \quad  \text{for Chebyshev-Gauss points}. \\
\end{cases}
\end{equation}

The spectral-Galerkin method makes use of a three-dimensional scalar product in the weighted $L^2_{\sigma}(\Omega)$ space, that is defined as
\begin{align}
\left<u, \upsilon \right>_{\sigma} &= \int_{\Omega} {u(\bm{x}) \upsilon^*(\bm{x})} \sigma(\bm{x})\,dxdydz, 
\end{align}
where $\upsilon^*$ is the complex conjugate of the test function $\upsilon$ and the weights $\sigma$ are unity for 
periodic directions and  $\sigma(x)=1/\sqrt{1-x^2}$ for the inhomogeneous direction. In 
this work we will make use of the discrete weighted $l^2_{\sigma}(\Omega)$ space, 
where quadrature is employed for the integration. As such, we redefine the scalar product as
\begin{equation}
\left<u, \upsilon\right>_{\sigma} = \sum_{i=0}^{N_x}\sum_{j=0}^{N_y-1} \sum_{k=0}^{N_z-1} u(x_i, y_j, z_k) \upsilon^*(x_i, y_j, z_k) \sigma(x_i), \label{eq:quadrature}
\end{equation}
that is more amendable to the fast transforms that will be defined later in this section.

The Navier Stokes equations are discretized using Fourier basis functions for the periodic directions, and a combination of Chebyshev polynomials in the wall-normal direction. Three different sets of basis functions and function spaces are relevant for the wall-normal direction
\begin{align}
&  \phi_k(x) = T_k(x), & W_{N_x} &= \text{span}\{\phi_k\}_{k=0}^{N_x}, \label{eq:Tk}\\
& \D{\phi}_k(x) = T_k(x) - T_{k+2}(x), & \D{W}_{N_x} &= \text{span} \{ \D{\phi}_k\}_{k=0}^{N_x-2}, \label{eq:phiD}\\
& \N{\phi}_k(x) = T_k(x) - \frac{2(k+2)}{k+3} T_{k+2}(x) + 
\frac{k+1}{k+3} T_{k+4}(x), & \N{W}_{N_x} &= \text{span} \{\N{\phi}_k\}_{k=0}^{N_x-4}, \label{eq:phiN} 
\end{align}
where $T_k(x)$ is the $k$'th degree Chebyshev polynomial of the first kind. The 
basis functions and function spaces in (\ref{eq:phiD}) and (\ref{eq:phiN}) were 
suggested by Shen \cite{Shen95}, and satisfy, respectively, the boundary conditions 
$\D{\phi}_k(\pm 1) = 0$, $\N{\phi}_k(\pm 1)=0$ and $\N{\phi}'_k(\pm 1)=0$. 

Three-dimensional basis functions and function spaces, that are periodic in $y$ 
and $z$ directions, can now be defined as
\begin{align}
  \psi_{\bm{k}}(\bm{x}) = \phi_{l}(x)e^{ \imath(\underline{m} y + \underline{n} z)}, \quad V_N &= \text{span} \{ \psi_{\bm{k}}(\bm{x}):\, \bm{k} \in K_N  \}, \\
  \D{\psi}_{\bm{k}}(\bm{x}) = \D{\phi}_{l}(x)e^{ \imath(\underline{m} y + \underline{n} z)}, \quad \D{V}_N &= \text{span} \{ \D{\psi}_{\bm{k}}(\bm{x}):\, \bm{k} \in \D{K}_N  \}, \\
  \N{\psi}_{\bm{k}}(\bm{x}) = \N{\phi}_{l}(x)e^{ \imath(\underline{m} y + \underline{n} z)}, \quad \N{V}_N &= \text{span} \{ \N{\psi}_{\bm{k}}(\bm{x}):\, \bm{k} \in \N{K}_N  \},
\end{align}
where $\imath=\sqrt{-1}$. 
The wavenumbermesh, $K_N$, for space $V_{N}$, is defined by the Cartesian product of wavenumbers from the two periodic and the inhomogeneous wall-normal direction: $K_N(l,m,n)= K^x(l) \times K^p(m,n)$, where
\begin{align}
K^p = \Big\{&(\underline{m}, \underline{n}) = \left(\frac{2 \pi m}{L_y}, \frac{2 \pi n}{L_z} \right), \text{where} \notag \\
&(m, n) \in  [-\frac{N_y}{2},-\frac{N_y}{2}+1,\ldots, 
\frac{N_y}{2}-1] \times [-\frac{N_z}{2},-\frac{N_z}{2},+1,\ldots, \frac{N_z}{2}-1] \Big\} 
\label{eq:wavenumbermeshp}
\end{align}
and $K^x(l) = \{l \in \mathbb{Z} | l=0,1, \ldots, N_x\}$. The two remaining wavenumber meshes, $\D{K}_{N}$ and $\N{K}_{N}$, differ from $K_N$ only in the range of the first index sets, $\D{K}^x$ and $\N{K}^x$, ending in $N_x-2$ and $N_x-4$, respectively (see Eqs. (\ref{eq:phiD}) and (\ref{eq:phiN})). 

In the spectral-Galerkin method we look for solutions of the velocity 
components of the form
\begin{align}
u(\bm{x}, t) &= \frac{1}{N_yN_z}\sum_{\bm{k} \in \N{K}_N} \hat{u}_{\bm{k}}(t) 
\N{\psi}_{\bm{k}}(\bm{x}), \label{eq:u_solx} \\
v(\bm{x}, t) &= \frac{1}{N_yN_z}\sum_{\bm{k} \in \D{K}_N} \hat{v}_{\bm{k}}(t) 
\D{\psi}_{\bm{k}}(\bm{x}), \label{eq:u_soly} \\
w(\bm{x}, t) &= \frac{1}{N_yN_z}\sum_{\bm{k} \in \D{K}_N} \hat{w}_{\bm{k}}(t) 
\D{\psi}_{\bm{k}}(\bm{x}), \label{eq:u_solz}
\end{align}
where $\hat{u}_{\bm{k}}(t) = \hat{u}(l, {m}, {n}, t)$ are the expansion 
coefficients for the velocity component in $x$-direction (and similar for the 
other two components) and the scaling by $N_y$ and $N_z$ is merely for 
convenience 
and compliance with the definition used later for the inverse discrete Fourier 
transform. Note that from now on we will simply use the notation $\hat{u}$ for 
$\hat{u}_{\bm{k}}(t)$, when it is possible to simplify without loss of clarity. 
Likewise we will simply use $u$ for $u(\bm{x}, t)$. 

For an efficient method it is crucial to be able to compute $u$ quickly from $\hat{u}$, or, vice versa, to compute $\hat{u}$ quickly from $u$. This may be achieved using the fast transform methods to be defined next.
Consider first how to compute $u$ from the known expansion coefficients $\hat{u}$. Writing out the summation terms, the expression (\ref{eq:u_solx}) may be 
evaluated on the entire mesh (\ref{eq:Xn}) as
\begin{align}
u(x_i, y_j, z_k, t) &= \underbrace{\frac{1}{N_z}\sum_n 
\underbrace{\frac{1}{N_y} \sum_m \underbrace{\sum_l 
\hat{u}(l,m,n,t)\N{\phi}_l(x_i)}_{\N{\mathcal{S}}_x^{-1}} e^{\imath 
\underline{m} 
y_j}}_{\mathcal{F}_y^{-1}} e^{\imath \underline{n} z_k}}_{\mathcal{F}_z^{-1}} 
\quad \forall \, \bm{x} \in X_N, \notag \\
  u &= \N{\mathcal{T}}^{-1}(\hat{u}) =  \mathcal{F}^{-1}_z (\mathcal{F}^{-1}_y 
  (\N{\mathcal{S}}^{-1}_x (\hat{u}))),
  \label{eq:ifft} 
\end{align}
where $\N{\mathcal{T}}^{-1}$ is used as short notation for the complete inverse 
transform in space $\N{V}_N$, and $\mathcal{F}_{y}^{-1}$ and $\mathcal{F}_{z}^{-1}$ represent inverse discrete
Fourier transforms along directions $y$ and $z$ respectively. For simplicity, we have introduce here a special notation called the inverse 
\emph{Shen} transform, ${\mathcal{S}}^{-1}_x$, that is used to transform coefficients from spectral to physical space in a series expansion that is using either one of the bases in (\ref{eq:Tk}, \ref{eq:phiD}, \ref{eq:phiN}). The inverse Shen transform, $\mathcal{S}^{-1}_x$, is performed along the wall-normal 
$x$ direction, and it may be computed using fast 
Chebyshev transforms for all the three bases in (\ref{eq:Tk}, \ref{eq:phiD}, 
\ref{eq:phiN}), see Alg~\ref{alg:ifst} in the Appendix. The transforms are slightly different for the three 
bases and 
${\mathcal{S}}_x^{-1}, \D{\mathcal{S}}_x^{-1}$ and $\N{\mathcal{S}}_x^{-1}$ are used 
to distinguish between them with obvious notation. Similarly, $\D{\mathcal{T}}^{-1}$ and ${\mathcal{T}}^{-1}$ define the inverse transforms for spaces $\D{V}_N$ and $V_N$. Note that a transform in  any one direction is performed over all indices in the other two directions. For example, for the transform in the $x$-direction we have
\begin{equation}
u(x_i,m,n,t) = \N{\mathcal{S}}_x^{-1}(\hat{u}(l,m,n,t)) \quad \forall \, i \in [0,1,\ldots, N_x] \text{ and }  m,n \in {K}^p,
\end{equation}
and similar for the other two directions.

Fast transforms may also be used in the scalar product defined in Eq. (\ref{eq:quadrature}), using basis $\N{\psi}_{\bm{k}}$ for $\upsilon$
\begin{align}
 \left< u, \N{\psi}_{\bm{k}} \right>_{\sigma} &= h \underbrace{\sum_i 
 \underbrace{\sum_j \underbrace{\sum_k u(x_i, y_j, z_k, t)  e^{-\imath 
 \underline{n} z_k}}_{\mathcal{F}_n}  e^{-\imath \underline{m} y_j} 
 }_{\mathcal{F}_m} \N{\phi}_{l}(x_i) \sigma(x_i)}_{\N{\mathcal{S}}_l} \quad \forall \, \bm{k} \in \N{K}_N,   \notag \\
  &=  h \N{\mathcal{S}}(u) = h \N{\mathcal{S}}_{l} (\mathcal{F}_{{m}} 
  (\mathcal{F}_{{n}}(u))). \label{eq:sst1}
\end{align}
Here $h\N{\mathcal{S}}(\cdot)$ denotes the complete three-dimensional scalar product 
and $h = L_yL_zN_y^{-1}N_z^{-1}$ is a constant. $\mathcal{F}_{{n}}$ and 
$\mathcal{F}_{{m}}$ represent discrete Fourier transforms in $z$- and 
$y$-directions, respectively, and $\N{\mathcal{S}}_{l}(\cdot) = (\cdot, 
\N{\phi}_l)_{\sigma}$ is used to represent the forward \emph{Shen} scalar product in the 
$x$-direction. The weights, $\sigma(x_i)$, required for the Shen transforms, are defined as
\begin{equation}
 \sigma(x_i) = \begin{cases}
       \frac{\pi}{c_i N_x} &\forall \, i=0,1,\ldots, N_x \quad  \text{for 
       Chebyshev-Gauss-Lobatto points},\\
       \frac{\pi}{N_x+1} &\forall \, i=0,1,\ldots, N_x  \quad \text{for 
       Chebyshev-Gauss points},      
 \end{cases}
\end{equation}
where $c_0 = c_{N_x} = 2$ and $c_i = 1$ for $0 < i < N_x$, see Sec. 1.11 of \cite{kopriva09}.

The scalar product  $h\N{\mathcal{S}}$ in (\ref{eq:sst1}) does not represent a 
complete transform. To find a transformation from physical $u$ to spectral 
$\hat{u}$, we make use of Eq. (\ref{eq:u_solx}) directly on the left hand side of (\ref{eq:sst1}) and use extensively the discrete 
orthogonality of the Fourier basis functions
\begin{align}
\left<u, \N{\psi}_{\bm{k}}\right>_{\sigma} &= \frac{1}{N_yN_z}\left< 
\sum_{(q,r,s) \in \N{K}_N} \hat{u}(q,r,s,t) \N{\psi}(q,r,s), 
\N{\psi}(l,m,n) \right>_{\sigma} \quad \forall \, \bm{k} \in \N{K}_N, \notag \\
           &= h \sum_{q=0}^{N_x-4} \sum_{i=0}^{N_x} \N{\phi}_{q}(x_i) 
           \N{\phi}_{l}(x_i) \sigma(x_i) \, \hat{u}(q, {m}, {n}, t), \notag \\
           &= h \sum_{q=0}^{N_x-4} \N{B}_{lq} \hat{u}(q, {m}, {n}, t). 
           \label{eq:sst2}
\end{align}
Here $\N{B}_{lq} = (\N{\phi}_q, \N{\phi}_l)_{\sigma} = 
\sum_{i=0}^{N_x} \N{\phi}_{q}(x_i) \N{\phi}_{l}(x_i) \sigma(x_i)$ are the components of a 
banded mass matrix with only 5 nonzero diagonals, see Table~\ref{tab:matrices}. The complete transformation is now obtained by 
setting Eq. (\ref{eq:sst2}) equal to (\ref{eq:sst1}) and solving for $\hat{u}$. 
Moving to matrix notation we get
\begin{align}
h\sum_{q=0}^{N_x-4} \N{B}_{lq}\, \hat{u}(q, {m}, {n}, t) &= 
h\N{\mathcal{S}}(u)(l, m, n, t)  \quad \forall \, \bm{k} \in \N{K}_N, \notag \\
 \hat{u} &= \N{\mathcal{T}}(u) =  \N{B}^{-1}\N{\mathcal{S}}(u), \label{eq:fstB}
\end{align}
where $\N{\mathcal{T}}(u)$ denotes the complete transformation, such that $u = 
\N{\mathcal{T}}^{-1}(\N{\mathcal{T}}(u))$. Note that, since $\N{B}$ assembles to 
a pentadiagonal matrix for its decoupled odd and even coefficients, the solution ($\N{B}^{-1}$) can be obtained very fast and the complete 
transformation thus requires a fast Chebyshev transform ($\mathcal{O}(N_x \log N_x)$) and 
a fast linear algebraic solve ($\mathcal{O}(N_x)$) for the wall-normal direction, given $m$ and $n$. Similar transforms $\mathcal{T}$ 
and $\D{\mathcal{T}}$ are defined for the two other bases (\ref{eq:Tk}) and 
(\ref{eq:phiD}), using mass matrices with components $B_{lq}=(T_q, T_l)_{\sigma}$ (diagonal) and 
$\D{B}_{lq}=(\D{\phi}_l, \D{\phi}_q)_{\sigma}$ (three non-zero diagonals) and scalar products 
$\mathcal{S}_l(\cdot) = (\cdot, T_l)_{\sigma}$ and $\D{\mathcal{S}}_l(\cdot) = (\cdot, 
\D{\phi}_l)_{\sigma}$. See Alg.~\ref{alg:fst} and \ref{alg:ifst} in the Appendix for algorithms for all required transforms.

\section{Discretization of Navier Stokes equations}
\label{sec:discretizationNS}
The Navier Stokes equations (\ref{eq:NS}) are solved using a scheme 
proposed by Kim, Moin and Moser \cite{Kim87}. This scheme is developed by taking the 
Laplacian of the wall-normal momentum equation and the curl of the momentum 
equation. Following elimination of the pressure, the equations to solve are 
\begin{align}
\frac{\partial}{\partial t} \nabla^2 u &= h_u + \nu \nabla^4 u, 
\label{eq:biharmonic} \\
\frac{\partial g}{\partial t} &= h_g + \nu \nabla^2 g, \label{eq:g} \\
f + \frac{\partial u}{\partial x} &= 0, \label{eq:f}
\end{align}
where
\begin{align}
f &= \frac{\partial v}{\partial y} + \frac{\partial w}{\partial z}, \\
g &= \frac{\partial w}{\partial y} - \frac{\partial v}{\partial z}, \\
h_u = -\frac{\partial}{\partial x} &\left( \frac{\partial 
\mathcal{H}_y}{\partial y} + \frac{\partial \mathcal{H}_z}{\partial z} \right) 
+ \left(\frac{\partial^2}{\partial y^2} + \frac{\partial^2}{\partial z^2} 
\right) \mathcal{H}_x ,
\\
h_g &= \frac{\partial \mathcal{H}_z}{\partial y} - \frac{\partial 
\mathcal{H}_y}{\partial z}.
\end{align}
The two remaining velocity components are computed from the definitions of $f$ 
and $g$. The biharmonic equation (\ref{eq:biharmonic}) is solved with boundary 
conditions $u(\pm 1) = u'(\pm 1) = 0$, where the Neumann conditions follow 
from the continuity equation. The boundary conditions for Eq. (\ref{eq:g}) are 
$g(\pm 1) = 0$. 

We consider the spectral-Galerkin discretization of (\ref{eq:biharmonic}), (\ref{eq:g}) and (\ref{eq:f}), using a central difference in time, with Crank-Nicolson for the linear terms and an Adams-Bashforth method for the nonlinear terms. To this end, the discretization in time is performed with a constant time step $\triangle t>0$, 
such that time is represented discretely as $t_{\kappa} = {\kappa} \triangle t, \, {\kappa} \in [0, 1, \ldots]$, and variables with a time superscript, like $u^{\kappa} = u(t_{\kappa})$. We get the variational formulation: with $u^0(\bm{x}), u^1(\bm{x}), g^0(\bm{x})$ and $g^1(\bm{x})$ given, for all $\kappa > 1$ find ${u}^{\kappa+1} \in \N{V}_N$, ${g}^{\kappa+1} \in \D{V}_N$ and ${f}^{\kappa+1} \in \D{V}_N$ such that
\begin{align}
  \left<  \frac{\nabla^2(u^{\kappa+1}-u^{\kappa})}{\triangle t}, \N{\psi}_{\bm{k}}\right>_{\sigma} &= 
\left<h_u^{\kappa+\half}, \N{\psi}_{\bm{k}} \right>_{\sigma} + {\nu} \left<\nabla^4 u^{\kappa+\half}, \N{\psi}_{\bm{k}}\right>_{\sigma} 
&\forall \,\psi_{\bm{k}} \in \N{V}_N, \label{eq:u1} \\
  \left<\frac{g^{\kappa+1}- g^{\kappa}}{\triangle t}, \D{\psi}_{\bm{k}}\right>_{\sigma} &= \left<h_g^{\kappa + \half}, 
\D{\psi}_{\bm{k}}\right>_{\sigma} + {\nu} 
\left<\nabla^2 g^{\kappa+\half}, \D{\psi}_{\bm{k}}\right>_{\sigma} &\forall \, \psi_{\bm{k}} \in \D{V}_N ,
\label{eq:g1} \\
	\left<f^{\kappa+1}, \D{\psi}_{\bm{k}}\right>_{\sigma} &= \left<\frac{\partial u^{\kappa+1}}{\partial x}, 
	\D{\psi}_{\bm{k}}\right>_{\sigma} &\forall \, \psi_{\bm{k}} \in \D{V}_N. \label{eq:f1}
\end{align}
Here superscript $\kappa + \half$ is used to represent Crank-Nicolson for linear terms (e.g., $u^{\kappa+\half}=0.5(u^{\kappa+1}+u^{\kappa})$) and Adams-Bashforth for nonlinear (e.g., $h_u^{\kappa + \half} = 1.5h_u^{\kappa} - 0.5 h_u^{\kappa-1}$).
 
With ${f}^{\kappa+1}$ and ${g}^{\kappa+1}$ known, the two remaining velocity components are then computed by projection to the 
Dirichlet space $\D{V}_N$: Find 
${v}^{\kappa+1} \in \D{V}_N$ and $w^{\kappa+1} \in \D{V}_N$ such that
\begin{align}
\left<f^{\kappa+1}, \D{\psi}_{\bm{k}}\right>_{\sigma} &= \left<\frac{\partial v^{\kappa+1}}{\partial y} + 
\frac{\partial w^{\kappa+1}}{\partial z}, \D{\psi}_{\bm{k}}\right>_{\sigma} \, &\forall \, \D{\psi}_{\bm{k}} \in 
\D{V}_N, \label{eq:f2} \\
\left<g^{\kappa+1}, \D{\psi}_{\bm{k}}\right>_{\sigma} &= \left<\frac{\partial w^{\kappa+1}}{\partial y}  - 
\frac{\partial v^{\kappa+1}}{\partial z}, \D{\psi}_{\bm{k}}\right>_{\sigma} \, &\forall \, \D{\psi}_{\bm{k}} \in 
\D{V}_N, \label{eq:g2}
\end{align}
which simplifies considerably because all the derivatives are in periodic 
directions. Written in spectral space Eqs. (\ref{eq:f2}) and (\ref{eq:g2}) 
become simply the algebraic expressions
\begin{align}
\hat{f}_{\bm{k}}^{\kappa+1} &= \imath \underline{m}\, \hat{v}_{\bm{k}}^{\kappa+1} + \imath 
\underline{n}\, \hat{w}_{\bm{k}}^{\kappa+1} &\forall \, \bm{k} \in \D{K}_N^0, \label{eq:f3} \\
\hat{g}_{\bm{k}}^{\kappa+1} &= \imath \underline{m}\, \hat{w}_{\bm{k}}^{\kappa+1} - \imath 
\underline{n}\, \hat{v}_{\bm{k}}^{\kappa+1} & \forall \, \bm{k} \in \D{K}_N^0, \label{eq:g3}
\end{align}
where $\D{K}_N^0$ is used to denote that these equations can be solved for all 
wavenumbers except $m=n=0$. For $m=n=0$ we solve instead the momentum equations in $y$ and $z$ directions (see Eq. (\ref{eq:NS})):
\begin{align}
\left< \frac{v^{\kappa+1} - v^{\kappa}}{\triangle t}, \psi_{\bm{k}} \right>_{\sigma} &= \left< \mathcal{H}_{y}^{\kappa+\half}, \psi_{\bm{k}} \right>_{\sigma} + \nu \left< \nabla^2 v^{\kappa+\half},  \psi_{\bm{k}} \right>_{\sigma} -\left< \beta^{\kappa+1}, \psi_{\bm{k}} \right>_{\sigma} & \forall \, \bm{k} \in \D{K}^x \times K^p(0,0), \label{eq:v0} \\
\left< \frac{w^{\kappa+1} - w^{\kappa}}{\triangle t}, \psi_{\bm{k}} \right>_{\sigma} &= \left< \mathcal{H}_{z}^{\kappa+\half}, \psi_{\bm{k}} \right>_{\sigma} + \nu \left< \nabla^2 w^{\kappa+\half},  \psi_{\bm{k}} \right>_{\sigma}  & \forall \, \bm{k} \in \D{K}^x \times K^p(0,0) \label{eq:w0}.
\end{align}
Note that the regular pressure is eliminated since $m=n=0$, and that (\ref{eq:v0}) is the only place where the driving force, or the mean pressure gradient, $\beta$, enters the equations. Also, since $\beta$ is constant in space, the term $<\beta^{\kappa+1}, \D{\psi}_{\bm{k}} >_{\sigma}$ will only be non-zero for $l=m=n=0$. 

The final step of the method is to rewrite all equations on matrix form, using one-dimensional scalar products. Inserting the expansion (\ref{eq:u_solx}) for $u$ in $\N{V}_N$, and similar for $g, h_u$ and $h_g$ in $\D{V}_N$, the inner products required to solve Eqs. (\ref{eq:u1}-\ref{eq:f1}) are
\begin{align}
\left<\nabla^4u, \N{\psi}_{\bm{k}}\right>_{\sigma} &= h\left[ \left( 
\N{\phi}_q^{''''}, 
\N{\phi}_l\right)_{\sigma} -2(\underline{m}^2+\underline{n}^2) \left( \N{\phi}_q^{''}, 
\N{\phi}_l\right)_{\sigma} + (\underline{m}^2+\underline{n}^2)^2\left( \N{\phi}_q, 
\N{\phi}_l\right)_{\sigma}  \right] \hat{u}_q, \\
\left< \nabla^2 u, \N{\psi}_{\bm{k}}\right>_{\sigma} &= h\left[\left( \N{\phi}_q^{''}, 
\N{\phi}_l\right)_{\sigma} - (\underline{m}^2+\underline{n}^2)\left( \N{\phi}_q, 
\N{\phi}_l \right)_{\sigma} \right] \hat{u}_q, \\
\left< \nabla^2 g, \D{\psi}_{\bm{k}}\right>_{\sigma} &= h\left[\left( \D{\phi}_q^{''}, 
\D{\phi}_l\right)_{\sigma} - (\underline{m}^2+\underline{n}^2)\left( \D{\phi}_q, 
\D{\phi}_l \right)_{\sigma} \right] \hat{g}_q, \\
\left<\frac{\partial u}{\partial x}, \D{\psi}_{\bm{k}}\right>_{\sigma} &=
h\left(\N{\phi}_q^{'}, \D{\phi}_l\right)_{\sigma} \hat{u}_q, \\
\left<h_u, \N{\psi}_{\bm{k}} \right>_{\sigma} &= 
h\left[-(\underline{m}^2+\underline{n}^2) 
\left(\D{\phi}_q, 
\N{\phi}_l \right)_{\sigma} \hat{\mathcal{H}}_{x, q} - \imath \underline{m}\left(\D{\phi}_q^{'}, 
\N{\phi}_l \right)_{\sigma} \hat{\mathcal{H}}_{y, q} - \imath 
\underline{n}\left(\D{\phi}_q^{'}, \N{\phi}_l \right)_{\sigma} 
\hat{\mathcal{H}}_{z, q}\right], \label{eq:S_hv}
\\
\left< h_g, \D{\psi}_{\bm{k}} \right>_{\sigma} &= h\left[ \imath \underline{m}\, \left(\D{\phi}_q, \D{\phi}_l \right)_{\sigma} 
\hat{\mathcal{H}}_{z, q} - \imath \underline{n}\, \left(\D{\phi}_q, \D{\phi}_l \right)_{\sigma} \hat{\mathcal{H}}_{y, q} 
\right],\label{eq:S_hg}
\end{align}
where for brevity in notation (as before) it is understood that the scalar products
act along the first dimension of the transformed variables, i.e.,  $ 
(\N{\phi}_q^{'}, \D{\phi}_l)_{\sigma} \hat{u}_q$ is short for the matrix 
vector product $ \sum_{q=0}^{N_x-4}( \N{\phi}_q^{'}, \D{\phi}_l)_{\sigma} 
\hat{u}(q, {m}, {n}, t)$, for all $m$ and $n$. The scalar products are used to set up linear systems 
of equations for the inhomogeneous wall-normal direction. All scalar products 
$(\cdot, \cdot)_{\sigma}$ can be represented by sparse matrices.
The required matrices with components $ \N{B}_{lq} = ( \N{\phi}_q, 
\N{\phi}_l)_{\sigma}, \D{B}_{lq} = ( \D{\phi}_q, 
\D{\phi}_l)_{\sigma}, \N{M}_{lq} = ( \D{\phi}_q, 
\N{\phi}_l)_{\sigma}, \D{C}_{lq} = (\N{\phi}_q^{'}, \D{\phi}_l)_{\sigma}$ and 
$ \N{C}_{lq} = (\D{\phi}_q^{'}, \N{\phi}_l)_{\sigma}$ are given in Table~\ref{tab:matrices}, $\N{A}_{lq}= -( \N{\phi}_q^{''}, 
\N{\phi}_l)_{\sigma}, \D{A}_{lq}= -( \D{\phi}_q^{''}, \D{\phi}_l)_{\sigma}$ are given in \cite{Shen95}, whereas $\N{Q}_{lq} = (\N{\phi}^{''''}_q, \N{\phi}_l)_{\sigma}$ is given below in Eqs (\ref{eq:N_Q})-(\ref{eq:rk}). Note that the matrices are computed using quadrature, which has some implications for Chebyshev-Gauss-Lobatto (GL) points, where the rows of the highest modes differ from those presented in Lemma 2.1 and 3.1 of \cite{Shen95}. This disagreement, that follows from inexact quadrature at the highest mode, explains the inclusion of the $c_{k+4}$ term for matrix $\N{B}$ and the $c_{k+2}$ term for matrix $\D{B}$. The matrix with components $B_{lq}=(T_q, T_l)_{\sigma}$ is also different in the highest mode from Eq. 2.7 of \cite{Shen95}, but agrees with Eqs. 1.135 and 1.136 of \cite{kopriva09}.

Assembling all scalar products, the final matrix form of Eqs. (\ref{eq:u1}), (\ref{eq:g1}) and (\ref{eq:f1}), that can be used to solve for $\hat{u}^{{\kappa}+1}$, $\hat{g}^{{\kappa}+1}$ and $\hat{f}^{\kappa+1}$,
given wavenumbers $m$ and $n$, are now found as
\begin{align}
\N{H}\hat{u}^{{\kappa}+1} & = \left(2\N{A} - 2\underline{z}^2\N{B} - \N{H} 
\right)\hat{u}^{{\kappa}} + \triangle t \N{\mathcal{S}}(h_v^{{\kappa}+1/2})  & 
\forall \, {\bm{k}} \in \N{K}_N, 
\label{eq:ufin}\\ 
\D{H} \hat{g}^{{\kappa}+1} &= \left(2 \D{B}-\D{H}\right) 
\hat{g}^{{\kappa}} + \triangle t \D{\mathcal{S}}(h_g^{{\kappa}+1/2})  &\forall \, 
{\bm{k}} \in 
\D{K}_N, \label{eq:gfin} \\
\D{B} \hat{f}^{\kappa+1} &= \D{C}\hat{u} &\forall \, \bm{k} \in \D{K_N}. \label{eq:ffin}
\end{align}
The coefficient matrices are given as
\begin{align}
\N{H}(m, n) &= -\frac{\nu \triangle t}{2}\N{Q} + \left( 1 + \nu 
\triangle t 
\underline{z}^2 \right) \N{A} - \frac{2\underline{z}^2 + \nu \triangle t 
\underline{z}^4}{2} \N{B} &\forall \, m, n \in K^p, 
\label{eq:Biharmonic_matrix} \\
\D{H}(m, n) &= -\frac{\nu \triangle 
t}{2}\D{A} + (1 + \frac{\nu \triangle t \underline{z}^2}{2}) \D{B} 
&\forall \, m, n \in K^p, \label{eq:Helmholtz_matrix}
\end{align}
where $\underline{z}^2 = \underline{m}^2 + \underline{n}^2$. Equations 
(\ref{eq:v0}) and (\ref{eq:w0}) are also written on matrix form as
\begin{align}
\D{H}\hat{v}^{{\kappa}+1} &= 
\left(2\D{B} - \D{H} \right)\hat{v}^{{\kappa}} + \triangle t\D{B}\hat{\mathcal{H}}_{y}^{{\kappa}+\half} - \triangle t \D{\mathcal{S}}(\beta^{\kappa+1}) & \forall \, l \in \D{K}^x, m=n=0, 
\label{eq:v00} \\
\D{H}\hat{w}^{{\kappa}+1} &= 
\left(2\D{B} - \D{H} \right)\hat{w}^{{\kappa}} + \triangle t\D{B}\hat{\mathcal{H}}_{z}^{{\kappa}+\half} & \forall \, l \in \D{K}^x, m=n=0. 
\label{eq:w00}
\end{align}

Finally, the nonlinear
terms are computed with the recipes given in Eqs. (\ref{eq:S_hv}, \ref{eq:S_hg}). To this end the required vector $\bm{\hat{\mathcal{H}}}$ is found by projecting to the Dirichlet vector space $\D{V}_N^3$, which corresponds to transforming the nonlinear cross product, evaluated in real space, over the entire mesh
\begin{equation}
{\hat{\mathcal{H}}}_i = \mathcal{\D{T}}((\bm{u} \times \bm{\omega})_i) \quad \forall \, i=x,y,z. \label{eq:H_hat}
\end{equation}
The curl, $\bm{\omega} = (g, \partial_z u - \partial_x w, \partial_y u - 
\partial_x v)$, is computed by projecting each individual term to its appropriate spectral space, before transforming back to physical space. 
To this end, the two terms $\partial_x v$ and $\partial_x w$ are projected to 
$V_N$ (requires $C_{lq} = (\D{\phi}_j^{'}, T_k)_{\sigma}$, see Tab. \ref{tab:matrices}), whereas the remaining $\partial_y u$ and $\partial_z u$ are projected to 
$\N{V}_N$. With compact notation using the transforms, we obtain
\begin{align}
\partial_x v(\bm{x}) &= \mathcal{T}^{-1}(B^{-1}C\hat{v}) &\forall \, \bm{x} \in X_N, \\
\partial_x w(\bm{x}) &= \mathcal{T}^{-1}(B^{-1}C\hat{w}) &\forall \, \bm{x} \in X_N, \\
\partial_y u(\bm{x}) &= \N{\mathcal{T}}^{-1}(\imath \underline{m}\hat{u}) &\forall \, \bm{x} \in X_N, \\
\partial_z u(\bm{x}) &= \N{\mathcal{T}}^{-1}(\imath \underline{n}\hat{u}) &\forall \, \bm{x} \in X_N.
\end{align}
Here, with a slight abuse of notation, the left hand side is simply representing the respective expressions evaluated on the quadrature points in the real mesh $X_N$. 
Note that the nonlinear term is also generally in need of 
dealiasing, at least for the two periodic directions, but this is not discussed further in this paper.

\section{Implementation}
\label{sec:implementation}
An outline of the solution procedure, used for the numerical method described 
in Sec.~\ref{sec:discretizationNS}, is given in Algorithm 
\ref{alg:NS}.\footnote{Note that throughout this paper we are using the 
pseudocode conventions of Kopriva \cite{kopriva09}, with some minor differences. 
A vectorization statement like $\{w_k\}_{k=0}^{N} \gets \{u_k\}_{k=0}^N$ indicates that components $u_k$, for $k=0,1,\ldots,N$, are copied from 
$u_k$ to $w_k$. Similar conventions apply to matrices, e.g., 
$\{U_{1,k}\}_{k=0}^N = \{V_{2,k}\}_{k=0}^N$ can be used for copying row 2 of 
$V$ to row 1 of $U$. Broadcasting is also implied, here meaning that for $\{w\}_{k=0}^{N} \gets c$, where $c$ is a scalar, all elements of $w_k$ from $k=0$ to $N$ gets the scalar value $c$.} The major computational cost in Alg. \ref{alg:NS} comes from  computing the 
nonlinear convection and solving for Eqs.~(\ref{eq:ufin}, \ref{eq:gfin}). We consider in this section the linear solvers.

\begin{algorithm}
	\caption{Solution procedure for Navier Stokes equations.    
	 \label{alg:NS}}
 	\DontPrintSemicolon	
    	Initialize $\bm{u}^0(\bm{x})$ and $\bm{u}^1(\bm{x})$ \;
    	Compute $\hat{\bm{u}}^0_{\bm{k}},  \hat{\bm{u}}^1_{\bm{k}}, \hat{g}^0_{\bm{k}}$ and $\hat{g}^1_{\bm{k}}$ \;
        Set parameters (e.g., mesh and viscosity) and end time $T$ \;
		Compute LU decompositions of $\N{H}(m,n), \D{H}(m,n) 
		\,\forall \, m, n \in K^p$ \;
        $t \gets \triangle t$\; 
		${\kappa} \gets 1$ \;
    	Compute nonlinear convection 
    	$\bm{\hat{\mathcal{H}}}_{\bm{k}}^{0}(\hat{\bm{u}}^0_{\bm{k}}, 
    	\hat{g}^0_{\bm{k}})$ \;
    	\While{$t<T$}{
	    	Compute nonlinear convection $\bm{\hat{\mathcal{H}}}_{\bm{k}}^{\kappa}$ \;
	    	{Compute right hand sides of Eqs.~}(\ref{eq:ufin}), 
	    	(\ref{eq:gfin}) \;
		    {Solve Eq.}(\ref{eq:ufin}) for 
		    $\hat{u}^{{\kappa}+1}_{\bm{k}} \, \forall \,\bm{k} \in \N{K}_N$ \;
		    {Solve Eqs.}(\ref{eq:gfin}) for 
		    $\hat{g}^{{\kappa}+1}_{\bm{k}} \, \forall\, \bm{k} \in \D{K}_N$ \;
		    {Solve Eqs.}(\ref{eq:ffin}) for 
		    $\hat{f}^{{\kappa}+1}_{\bm{k}} \, \forall\, \bm{k} \in \D{K}_N$ \;
		    {Solve Eqs.}(\ref{eq:f3}, \ref{eq:g3}) for  
		    $\hat{v}_{\bm{k}}, \hat{w}_{\bm{k}}\, \forall\, \bm{k} 
		    \in \D{K}_N^0$ \;
		    {Solve Eqs.}(\ref{eq:v00}, \ref{eq:w00}) for 
		    $\hat{v}^{{\kappa}+1}, \hat{w}^{{\kappa}+1}\, 
		    \forall \, l \in \D{K}^x, m=n=0$ \;
		    ${\kappa} \gets {\kappa} + 1$ \;
		    $t \gets t+\triangle t$ \;
		    {Update to new time step} \;
	    }
\end{algorithm}

Shen \cite{Shen95} writes that it is possible to solve Eq.~(\ref{eq:gfin}) directly 
with essentially the same number of operations as a pentadiagonal solver. To this end, note that the matrix 
$\D{H}$ decouples into odd and even components, leading to two matrices of type upper Hessenberg. A direct LU decomposition (see, e.g., \cite{stewart98}), leads to a lower triangular  
matrix $\D{L}$ with only one subdiagonal and an upper triangular matrix $\D{U}$ that is dense. However, each row in $\D{U}$ 
contains at most three distinct values at $\D{U}_{k,k}, \D{U}_{k,k+2}$ and $\D{U}_{k,k+4}$, 
and then $\D{U}_{k,j} = \D{U}_{k,k+4}\, \forall \, j = k+6, k+8, \ldots N_x-2$. 
Consequently 
only three diagonals in $\D{U}$ need storage and the back solve can be traversed 
very efficiently in $\mathcal{O}(N_x)$. Note also that the decoupling into 
odd and even coefficients leads to two subsystem that may be trivially 
solved simultaneously in two threads. For optimal 
performance, the odd and even coefficient would then need to be stored 
contiguously in computer memory, and not alternately, which would be the normal way of storing the coefficients.

The biharmonic problem in Eq.~(\ref{eq:ufin}) is more challenging to solve for efficiently, but 
it is still possible, as suggested yet not devised by Shen \cite{Shen95}, to find an algorithm that is $\mathcal{O}(N_x)$ for a system of $N_x$ unknowns. Note 
that $\N{H}$ is the sum of three matrices $\N{Q}, \N{A}$ and $\N{B}$, and can be 
written as
\begin{equation}
\N{H} = \xi_0\N{Q} + \xi_1\N{A} + \xi_2 \N{B},
\end{equation}
where $\xi_0, \xi_1$ and $\xi_2$ are constants (depending on $m$ and $n$). The matrix $\{\N{A}_{kj}\}_{k,j=0}^{N_x-4}$ contains only 
three nonzero diagonals at $ j = k-2, k, k+2$,  whereas matrix $\{\N{B}_{kj}\}_{k,j=0}^{N_x-4}$ 
contains five 
nonzero diagonals at $ j = k-4, k-2, k, k+2, k+4$. The nonzero elements of the 
upper triangular $\{\N{Q}_{kj}\}_{k,j=0}^{N_x-4}$ matrix are given as \cite{Shen95}
\begin{align}
\label{eq:N_Q}
 \N{Q}_{kk} &= -4(k+1)(k+2)^2, \\
 \N{Q}_{kj} &= p_kq_j + r_ks_j, \quad j = k+2, k+4, \ldots, N_x-4,
\end{align}
where 
\begin{align}
p_k &= 8 k (k+1)(k+2)(k+4)\pi, &q_j &= \frac{1}{j+3}, \label{eq:pk} \\
r_k &= 24(k+1)(k+2)\pi, &s_j &=  \frac{(j+2)^2}{(j+3)}. \label{eq:rk}
\end{align} 
A straight forward direct LU decomposition (without pivoting) of $\N{H}$ can be performed as shown in Alg. 
\ref{alg:lu_biharmonic}, which leads to a lower 
triangular 
matrix $\{\N{L}_{kj}\}_{k,j=0}^{N_x-4}$ with two nonzero diagonals at $j=k-2$ and $k-4$ plus a unity 
main diagonal. The upper triangular matrix $\N{U}$ is dense and as such a 
show stopper for an efficient back solve. However, we note that $\{\N{U}_{kj}\}_{k,j=0}^{N_x-4}$ 
contains three distinct diagonals at $j=k, k+2$ and $k+4$, whereas the remaining part can be expressed as
\begin{equation}
\N{U}_{kj} = \xi_0(a_k q_j + b_k s_j), \quad j = k+6, k+8, \ldots, N_x-4 \text{ and } k\le N_x-10, 
\label{eq:ab}
\end{equation}
where $\{a_k\}_{k=0}^{N_x-10}$ and $\{b_k\}_{k=0}^{N_x-10}$ are two new 
vectors that can be computed recursively from $\N{L}$, as shown in Alg.~\ref{alg:ab}. 
If $a$ and $b$
are pre-computed, 
the total storage requirement for the complete LU decomposition is less than 
$7N_x$, since there are two diagonals for $\N{L}$, three for $\N{U}$ plus ${a}$ and 
${b}$. The solution of the complete system can be obtained very quickly with 
one simple forward elimination and a back solve, where the back 
solve is very fast ($\mathcal{O}(N_x)$) because of (\ref{eq:ab}), leading 
to a formula for the backwards substitution of $\N{U}_{kj} \hat{u}_j = y_k$ like (valid for $k \le N_x-10$)
\begin{equation}
\hat{u}_k = \left(y_k - \N{U}_{k, k+2} \hat{u}_{k+2} - \N{U}_{k, k+4} \hat{u}_{k+4} - 
\xi_0 a_k\sum_{\substack{j=k+6\\k-j \text{ even}}}^{N_x-4} q_j \hat{u}_j - 
\xi_0 b_k\sum_{\substack{j=k+6\\k-j \text{ even}}}^{N_x-4} 
s_j \hat{u}_j\right)/\N{U}_{kk},
\end{equation}
that only requires 
one new addition to the row-sums for each back-traversed
row, see Alg. (\ref{alg:SolveBiharmonic}).

\SetKwFunction{LUbiharmonic}{LUbiharmonic}
\begin{algorithm}
    \caption{LU decomposition of biharmonic operator $\N{H}$. Current algorithm 
        is using dense storage of $\N{U}_{kj}$, but no more than the 3 items 
        per row need 
        to be stored ($\N{U}_{kk}, \N{U}_{kk+2}$ and $\N{U}_{kk+4}$) after the row is no longer needed in the for-loop below 
        (as the loop is traversed to row $k$, $\N{U}_{m, j}$ is no longer 
        needed for $m<k-4$). 
        \label{alg:lu_biharmonic}}
        \Fn{\LUbiharmonic{$\{\N{H}_{kj}\}_{k, j = 0}^{N}$}}{
        \Output{$\{\N{L}_{kj}\}_{k, j = 0}^{N}$, $\{\N{U}_{kj}\}_{k, j = 
        0}^{N}$}
        integer $M_e \gets \lfloor N/2 \rfloor$ \;
        integer $M_o \gets \lfloor (N-1)/2 \rfloor$  \;
        $\{\hat{U}_{0,2j}\}_{j=0}^{M_e} \gets \{\hat{H}_{0,2j}\}_{j=0}^{M_e}$ \;
        $\{\hat{U}_{1,2j+1}\}_{j=0}^{M_o} \gets 
        \{\hat{H}_{1,2j+1}\}_{j=0}^{M_o}$ \;
        $\N{L}_{20} \gets \N{H}_{20}/\N{U}_{00}$ \;
        $\N{L}_{31} \gets \N{H}_{31}/\N{U}_{11}$ \;
        $\{\hat{U}_{2,2j}\}_{j=1}^{M_e} \gets \{\N{H}_{2,2j}\}_{j=1}^{M_e} - 
        \N{L}_{20} \{\N{U}_{0,2j}\}_{j=1}^{M_e}$ \;
        $\{\hat{U}_{3,2j+1}\}_{j=1}^{M_o} \gets \{\N{H}_{3,2j+1}\}_{j=1}^{M_o} 
        - \N{L}_{31} \{\N{U}_{1,2j+1}\}_{j=1}^{M_o}$ \;
        
        \For{$k = 4$ \KwTo $N$}{
        integer $M \gets \lfloor \frac{N-k}{2}\rfloor$ \;
        $\N{L}_{k,k-4} \gets \N{H}_{k, k-4} / \N{U}_{k-4, k-4}$ \;
        $\{\N{U}_{k,k+2j}\}_{j=-1}^{M} \gets \{\N{H}_{k,k+2j}\}_{j=-1}^{M} - 
        \N{L}_{k, k-4} 
        \{\N{U}_{k-4, k+2j}\}_{j=-1}^{M}$ \;
        $\N{L}_{k,k-2} \gets {\N{U}}_{k,k-2}/{\N{U}}_{k-2, k-2}$ \;
        $\{\N{U}_{k,k+2j}\}_{j=-1}^{M} \gets \{{\N{U}}_{k,k+2j}\}_{j=-1}^{M} - 
        \N{L}_{k, k-2} \{\N{U}_{k-2,k+2j}\}_{j=-1}^{M}$ \;
        }
        \Return $\{\N{L}_{kj}\}_{k, j = 0}^{N}$, $\{\N{U}_{kj}\}_{k, j = 0}^{N}$
    }
\end{algorithm}	

\SetKwFunction{reduceAB}{ReduceAB}
\begin{algorithm}
	\DontPrintSemicolon
    \caption{Recursive formula to compute $\{a_k\}_{k=0}^{N_x-10}, 
	\{b_k\}_{k=0}^{N_x-10}$ in Eq. (\ref{eq:ab}). Note that the lower triangular 
	matrix $\{L_{kj}\}_{k,j=0}^{N_x-4}$ contains only two nonzero subdiagonals and 
	the parameters $p_k$ and $r_k$ are given in Eqs. (\ref{eq:pk}) and 
	(\ref{eq:rk}) respectively. \label{alg:ab}}
		\Fn{\reduceAB{$ \{\N{L}_{kj}\}_{k,j=0}^{N_x-4}$}}{
        \Output{$\{a_k\}_{k=0}^{N_x-10}, \{b_k\}_{k=0}^{N_x-10}$} 
        
        \For{$k=0$ \KwTo $1$}{
        	$a_k \gets p_k$ \;
        	$b_k \gets r_k$ \;
        	$a_{k+2} \gets p_{k+2} - \N{L}_{k+2,k}\,a_{k}$ \;
        	$b_{k+2} \gets r_{k+2} - \N{L}_{k+2,k}\,b_{k}$ \;
        }
        
		\For{$k=4$ \KwTo $N_x-10$}{
     		$a_k \gets p_k - \N{L}_{k,k-2}\,a_{k-2} - \N{L}_{k, k-4} \,a_{k-4}$ 
     		\;
	    	$b_k \gets r_k - \N{L}_{k, k-2}\,b_{k-2}- \N{L}_{k, k-4} \,b_{k-4}$ 
	    	\;
		}
		\Return $\{a_k\}_{k=0}^{N_x-10}, \{b_k\}_{k=0}^{N_x-10}$ \;
    }
\end{algorithm}

\SetKwFunction{LUSolveBiharmonic}{LUSolveBiharmonic}
\begin{algorithm}
    \DontPrintSemicolon
	\caption{Solve biharmonic Eq. (\ref{eq:ufin}) with pre-computed 
	$\{\N{U}_{kj}\}_{k,j=0}^{N}, 
		\{\N{L}_{kj}\}_{k,j=0}^{N}$ 
	matrices, as well as the $\{ a_k\}_{k=0}^{N-6}$ and $\{ b_k\}_{k=0}^{N-6}$ vectors. 
	The parameters $q_j$ and $s_j$ are given in Eqs. (\ref{eq:pk}) and 
	(\ref{eq:rk}) respectively. The vector $\{f_k\}_{k=0}^{N}$ contains the right hand side and $\{ \hat{u}_k\}_{k=0}^N$ is the solution.	\label{alg:SolveBiharmonic}}
		\Fn{\LUSolveBiharmonic{$\{f_k\}_{k=0}^{N}, \{\N{U}_{kj}\}_{k,j=0}^{N}, 
                \{\N{L}_{kj}\}_{k,j=0}^{N}, 
                \{a_k\}_{k=0}^{N-6}, \{b_k\}_{k=0}^{N-6}, 
                \xi_0$}}{
          \Output{$\{ \hat{u}_k\}_{k=0}^N$}
		  \tcp*[l]{Solve $\N{L}_{kj}y_j=f_k$ by forward elimination}
          \For{$k=0$ \KwTo $1$}
		  {$y_k \gets f_k$ \;
		   $y_{k+2} \gets f_{k+2} - \N{L}_{k+2,k} y_{k}$} 
		  \For{$k = 4$ \KwTo $N$}{
		  $y_k \gets f_k - \N{L}_{k,k-2}y_{k-2} - \N{L}_{k, k-4}y_{k-4}$ \;
		  }
		  \tcp*[l]{Solve $\N{U}_{kj}\hat{u}_j=y_k$ with back 
              substitution}      
          \For{$k = N$ \KwTo $N-1$ \KwStep $-1$}{
            $\hat{u}_k \gets y_k / \N{U}_{k,k}$ \;
            $\hat{u}_{k-2} \gets \left(y_{k-2} - \N{U}_{k-2, 
            	k}\hat{u}_k\right)/\N{U}_{k-2,k-2}$ \;
          }              
	      
	      $q^o \gets 0.0$ \;
          $q^e \gets 0.0$ \;
          $s^o \gets 0.0$ \;
          $s^e \gets 0.0$ \;
		  \For{$k = N-4$ \KwTo $0$ \KwStep $-1$}{
		  $\hat{u}_k \gets y_k - \N{U}_{k, k+2} \hat{u}_{k+2} - \N{U}_{k, k+4} 
		  \hat{u}_{k+4} $ \;
		  \If{$k < N-5$}{
		  $j \gets k+6$ \;
		  \If{$k$ is odd}{
		    $q^o \gets q^o + \hat{u}_{j}q_{j}$ \;
		    $s^o \gets s^o + \hat{u}_{j}s_{j}$ \;
		    $\hat{u}_k \gets \hat{u}_k - \xi_0q^o a_k - \xi_0s^o b_k $ \;
        }
		  \Else{
		    $q^e \gets q^e + \hat{u}_{j}q_{j}$ \;
		    $s^e \gets s^e + \hat{u}_{j}s_{j}$ \;	  
		    $\hat{u}_k \gets \hat{u}_k - \xi_0q^e a_k - \xi_0s^e b_k $ \;
		}
		}
		  $\hat{u}_k \gets  \hat{u}_k / \N{U}_{kk}$ \;
		  }
		  \Return $\{\hat{u}_k\}_{k=0}^{N}$ \;  
		}
\end{algorithm}

\section{Verification and validation}
\label{sec:verification}
We have in previous sections described a Navier-Stokes solver for channel flows applicable to large-scale simulations. We have devised algorithms for fast direct sparse solvers, scalar products, and for the necessary fast transforms between physical and spectral space. The complete solver has been implemented in the open source code spectralDNS~\cite{spectralDNS}. For MPI we are using a slab decomposition, as described by Mortensen \cite{Mortensen2016}, but parallel scalability is not considered here since all algorithms described in previous sections are executed in serial. In this section we will simply verify the implementation, and validate the method for large-scale simulations. To this end we will first study some simple 1D problems. 

Roundoff errors are often causing problems for spectral methods. The roundoff error of the linear solver described for the biharmonic problem in Sec.~\ref{sec:implementation}, and similarly for the Helmholtz problem, may be estimated 
by the following approach: Let $\{u_i\}_{i=0}^{N_x-4}$ be a uniformly distributed 
random vector of double precision in the interval $(0, 1)$. Compute the matrix vector product ${f}=\N{H}{u}$ and 
then solve $\N{H} {v} = {f}$, using Algorithms \ref{alg:lu_biharmonic}, \ref{alg:ab} and \ref{alg:SolveBiharmonic}, for 
$\{v_i\}_{i=0}^{N_x-4}$. The roundoff error of the solution 
algorithm may then be found as $\max \{|u_i-v_i|\}_{i=0}^{N_x-4} / \max \{|u_i|\}_{i=0}^{N_x-4}$. For this experiment we use double precision and assume that $\triangle t= 10^{-5}$ 
and $\nu=1/5200$, that are reasonable parameters for a turbulent channel flow 
simulation at ${Re}_{\tau}=5200$. Four different wavenumbers 
$\underline{z}=(0, 200, 1800, 5400)$ are 
also chosen as representative for such large-scale simulations. The 1800 case corresponds to the highest wavenumber used for the resolution 
of the largest known channel simulations, performed by Lee and 
Moser~\cite{leemoser15}, where maximum $\underline{z} = \sqrt{\underline{m}_{max}^2+\underline{n}_{max}^2} = \sqrt{(10240/8)^2 + (7680/6)^2} \approx 1800$. The mesh in the wall-normal direction 
is varied from 64 to 4096 points, which also covers almost three times higher 
resolution than used by Lee and Moser. Results shown in the 
first three rows of Table ~\ref{tab:roundoff} indicate that 
roundoff errors are not very pronounced, even for problems as large as 
$N_x=4096$. Further, the roundoff errors increase only slowly with increasing 
the wavenumbers. Corresponding results are also computed for the Helmholtz 
matrix (\ref{eq:Helmholtz_matrix}) and 
results are shown in the last three rows of Table \ref{tab:roundoff}. For 
the Helmholtz matrix roundoff errors are seen to be insignificant. Note that Chebyshev-Gauss points are used for the computations in Table~\ref{tab:roundoff} and that Chebyshev-Gauss-Lobatto points give very similar (slightly better) results.

\begin{table}
	\caption{ Roundoff error of the Biharmonic and Helmholtz linear algebra solvers. The first three 
	rows show results for Alg. (\ref{alg:SolveBiharmonic}) and matrix $\N{H}(\underline{z}, 
	\triangle t, \nu)$, whereas the last three rows show results for matrix 
	$\D{H}(\underline{z}, \triangle t, \nu)$. Parameters 
	are chosen as $\nu=1/5200$ and $\triangle t=10^{-5}$ and the numbers shown 
	are for an average of 100 runs.	 \label{tab:roundoff}}
\centering
   \begin{tabular}{cccccccc}
	   \multicolumn{8}{c}{$N_x$} \\
	   \hline
$\underline{z}$ & 64 & 128 & 256 & 512 & 1024 & 2048 & 4096 \\ 
\hline
0 & 3.27e-15  & 8.36e-15  & 2.10e-14  & 2.65e-14  & 2.88e-14  & 2.83e-14  & 3.29e-14  \\ 
200 & 5.77e-14  & 8.50e-14  & 1.03e-13  & 1.08e-13  & 1.06e-13  & 1.06e-13  & 1.14e-13  \\ 
1800 & 2.37e-13  & 1.65e-12  & 3.41e-12  & 3.33e-12  & 3.94e-12  & 3.79e-12  & 3.35e-12  \\ 
5400 & 1.76e-13  & 1.84e-12  & 1.24e-11  & 1.57e-11  & 1.58e-11  & 1.55e-11  & 1.51e-11  \\ 
\hline
0 & 3.26e-15  & 7.83e-15  & 1.89e-14  & 2.47e-14  & 2.16e-14  & 2.31e-14  & 1.94e-14  \\ 
200 & 3.34e-15  & 8.63e-15  & 2.04e-14  & 2.54e-14  & 2.19e-14  & 2.01e-14  & 2.08e-14  \\ 
1800 & 3.09e-15  & 8.23e-15  & 1.82e-14  & 2.27e-14  & 2.04e-14  & 2.18e-14  & 2.12e-14  \\ 
5400 & 3.15e-15  & 7.77e-15  & 1.88e-14  & 2.25e-14  & 2.22e-14  & 2.10e-14  & 2.02e-14 
		
\end{tabular}
\end{table}

To show that the solvers are scaling with size, we investigate the execution time 
for solving $\N{H} {v} = {f}$ and $\D{H} {v} = {f}$ on a complete 
wavenumber mesh of size $(N_x,64,64)$. We use a MacBook Pro laptop with the 2.8 GHz intel Core i7 processor. The average time for one single solve is shown in Table~\ref{tab:timings}, where the numbers in parenthesis are used to indicate scaling (unity would be perfectly $\mathcal{O}(N_x)$). The scaling is seen to be close to linear, and variations 
are most likely due to how well datastructures fit in the cache. The ability to solve large-scale problems is evident.
\begin{table}
	\centering
	\caption{Execution time (in seconds) for one solve  of the Helmholtz 
	(\ref{eq:gfin}) and Biharmonic 
	(\ref{eq:ufin}) problems. The numbers in parenthesis are computed as $t_n / t_{n-1} / 2$, where $t_n$ is the current solvers timing in row $n$. 	\label{tab:timings}}
	\begin{tabular}{ccc}
$N_x$ & Biharmonic & Helmholtz \\ 
\hline
64 & 5.15e-07 (0.00) & 4.01e-07 (0.00) \\ 
128 & 1.05e-06 (1.02) & 8.03e-07 (1.00) \\ 
256 & 2.13e-06 (1.01) & 1.65e-06 (1.03) \\ 
512 & 4.39e-06 (1.03) & 3.39e-06 (1.03) \\ 
1024 & 8.67e-06 (0.99) & 6.98e-06 (1.03) \\ 
2048 & 1.75e-05 (1.01) & 1.39e-05 (0.99) \\ 
4096 & 3.96e-05 (1.13) & 3.03e-05 (1.09) \\ 
8192 & 7.95e-05 (1.00) & 6.27e-05 (1.03) \\
	\end{tabular}
\end{table}

The Orr-Sommerfeld equation (see, e.g., \cite{Orszag_1971}) is further used to establish that we are solving the right equations, and that the accuracy of the solver is second-order in time. The Orr-Sommerfeld eigenvalue problem is first solved for $Re=8000$, using the biharmonic Shen basis (\ref{eq:phiN}). The leading eigenvalue/eigenvector pair is then added as a perturbation to a laminar, parabolic, base flow, leading to an initial field
\begin{align}
u(x,y,t) &= -\epsilon  \Re\{\imath \xi(x) \exp (\imath (y-\lambda t))\}, \label{eq:OS1} \\
v(x,y,t) &= 1-x^2 + \epsilon \Re\{\xi'(x) \exp(\imath(y-\lambda t)\}, \label{eq:OS2}
\end{align}
where $\xi(x)$ is the eigenvector, $\lambda=0.2470750602+\imath 2.664410371\cdot 10^{-3}$ the eigenvalue and $\epsilon \in \mathbb{R}^+$ is a small positive constant. Note that if $\epsilon$ is chosen as too big, then the initial/boundary value problem will not agree well with the linearized Navier-Stokes equations used in deriving the Orr-Sommerfeld equations, due to second-order terms. Likewise, if $\epsilon$ is chosen as too small, then roundoff errors will be too large, since the perturbation is added to a base flow with magnitude of order unity. It turns our that $\epsilon=10^{-7}$ is a good compromise, and that a mesh of size $128 \times 8 \times 2$ is accurate enough to isolate temporal errors. Since the evolution of the perturbation is known analytically, the $L^2$-norm of the deviation from this exact linearized solution is used as a measure of the error. We also compute the error in the perturbation energy by
\begin{equation}
E(t) =  \frac{\left<\tilde{\bm{u}}(\bm{x}, t), \tilde{\bm{u}}(\bm{x}, t) \right>}{\left<\tilde{\bm{u}}(\bm{x}, 0), \tilde{\bm{u}}(\bm{x}, 0) \right>}  - \exp(2 \Im\{\lambda\}t), \label{eq:OSenergy}
\end{equation}
where $\tilde{\bm{u}} = [u, v-(1-x^2)]$ represents the perturbation vector. 
Table \ref{tab:OS} shows the computed errors using a final time of $t=50$, corresponding to 2 flow-throughs, and a varying $\triangle t$. We observe, as expected, second-order accuracy in the $L^2$-norm in time. The error in the accumulated perturbation energy is also found to be close to second order. Table \ref{tab:OS_spatial} shows the computed error as a function of spatial discretization, where the time step is kept constant at $\triangle t=10^{-3}$ for 50 steps. The spectral decay of the error, with increasing $N_x$, is evident, and there is little difference between using either Chebyshev-Gauss (GC) or Chebyshev-Gauss-Lobatto (GL) points. Note that there is no difference at all between GC and GL for the results shown Table \ref{tab:OS}.

\begin{table}
	\centering
	\caption{Error in the Orr-Sommerfeld perturbation at $Re=8000$, for the most unstable eigenmode, as a function of time step. The mesh size is $128 \times 8 \times 2$. The $L^2$-norm $ \| \bar{\bm{u}} \| = \sqrt{\left< \bar{\bm{u}}, \bar{\bm{u}}\right>}$, where $\bar{\bm{u}}=\bm{u}-\bm{u}_{\mathrm{exact}}$, and the energy error is computed as Eq. (\ref{eq:OSenergy}). The numbers in parenthesis are showing the order of the error, computed as $\log(\tau_i/\tau_{i-1})/\log(\triangle t_i/ \triangle t_{i-1})$, where $i$ is the row number of the table and $\tau$ represents the values of $\|\bar{\bm{u}}\|$ or $E(t=50)$.  \label{tab:OS}}
	\begin{tabular}{ccc}	
 $\triangle t$ & $\| \bar{\bm{u}} \|(t=50)$ & $\int_{t=0}^{t=50}E(t) \mathrm{d}t$ \\
 \hline

1.0000e-01 & 2.9005e-09  (0.000) & 1.2775e-02 (0.000) \\
6.6667e-02 & 1.2877e-09  (2.003) & 5.3068e-03 (2.167) \\
5.0000e-02 & 7.2559e-10  (1.994) & 2.8939e-03 (2.108) \\
4.0000e-02 & 4.6445e-10  (1.999) & 1.8138e-03 (2.094) \\
3.3333e-02 & 3.2257e-10  (1.999) & 1.2418e-03 (2.078) \\
2.8571e-02 & 2.3701e-10  (1.999) & 9.0307e-04 (2.066) \\
2.5000e-02 & 1.8148e-10  (1.999) & 6.8608e-04 (2.058)
	\end{tabular}
\end{table}

\begin{table}
	\centering
	\caption{Error in the Orr-Sommerfeld perturbation at $Re=8000$, for the most unstable eigenmode, as a function of $N_x$. The time step $\triangle t = 10^{-3}$ and the end time is $t=0.05$. The $L^2$-norm $ \| \bar{\bm{u}} \| = \sqrt{\left< \bar{\bm{u}}, \bar{\bm{u}}\right>}$, where $\bar{\bm{u}}=\bm{u}-\bm{u}_{\mathrm{exact}}$, and the energy error is computed as Eq. (\ref{eq:OSenergy}). GC and GL refer to Chebyshev-Gauss and Chebyshev-Gauss-Lobatto points, respectively.   \label{tab:OS_spatial}}
	\begin{tabular}{ccccc}	
		$N_x$ & \multicolumn{2}{c}{$\frac{1}{\epsilon}{\| \bar{\bm{u}} \|(t=0.05)}$} & \multicolumn{2}{c}{$\int_{t=0}^{t=0.05}E(t) \mathrm{d}t$} \\
		\hline
			 & GC & GL & GC & GL \\
		\hline
16	&	3.23081791e-01	&	4.76365364e-01	&	3.35105681e-02	&	3.10940983e-03	\\
32	&	5.71635963e-03	&	6.77845098e-03	&	2.57328408e-04	&	5.02990613e-05	\\
64	&	6.99681587e-08	&	8.22331688e-08	&	1.75070691e-09	&	8.04535105e-10	\\
128	&	5.06389229e-08	&	4.81745061e-08	&	1.01924047e-09	&	1.53428292e-09	\\
256	&	4.89447535e-08	&	4.83040873e-08	&	8.40445491e-10	&	3.77475917e-09	
	\end{tabular}
\end{table}

The final test for this fully spectral Navier-Stokes solver is the pressure driven turbulent channel flow. First, to measure scaling of the complete 3-dimensional solver, we use a rather small turbulent Reynolds number of $Re_{\tau}=590$ and a small computational box of size $N_x \times 32 \times 32$.  Since the resolutions in $y$ and $z$ directions are kept constant, the complete Navier Stokes solver should ultimately scale like and 
$\mathcal{O}(N_x\log N_x)$, since it is a combination of linear algebra routines 
$\mathcal{O}(N_x)$, fast transforms $\mathcal{O}(N_x \log N_x)$ and several elementwise array operations $\mathcal{O}(N_x)$. We show in Table \ref{tab:timings_solver} some timings for one single time 
step, increasing only the wall-normal discretization. We see that the total 
time scales approximately as $N_x \log N_x$ for large $N_x$ and linearly 
for small $N_x$. This shift happens since at small $N_x$ the transforms are fast and linear operations, like elementwise multiplications, are more dominating. 
We also see that the time for assembly of the right hand side is approximately 
$\mathcal{O}(N_x \log N_x)$, whereas the time for solving the linear systems is 
very nearly linear. Note that for fast transforms we are using 
the FFTW library \cite{libfftw}, with planner effort set to FFTW$\_$MEASURE. 
The results in Table \ref{tab:timings_solver} are computed with GC points and dealiasing is
performed with the 2/3-rule \cite{Orzag71}. 

\begin{table}
	\centering
	\caption{Computing time for one time step using 32 wavenumbers in each periodic direction and $N_x$ points in the wall-normal direction. The Assemble and Solve columns show the time required to assemble the right hand side (should be $\mathcal{O}(N_x \log N_x)$ due to transforms) and solve the linear systems, respectively. The number in parenthesis shows the scaling, and unity would mean $\mathcal{O}(N_x \log N_x)$ for Total and Assemble, whereas it would indicate $\mathcal{O}(N_x)$ for Solve. \label{tab:timings_solver}}
	\begin{tabular}{cccc}	
		$N_x$ & Total ($\frac{t_k }{t_{k-1}} \frac{\log N_x-1}{2 \log N_x}$) & Assemble ($\frac{t_k }{t_{k-1}} \frac{\log N_x-1}{2 \log N_x}$) & Solve ($\frac{t_k}{t_{k-1}} \frac{1}{2}$) \\ 
		\hline
		32 & 0.010 (0.00) & 7.83e-03 (0.00) & 1.07e-03 (0.00) \\ 
		64 & 0.019 (0.86) & 1.53e-02 (0.88) & 1.98e-03 (0.92) \\ 
		128 & 0.040 (0.96) & 3.23e-02 (0.97) & 3.81e-03 (0.96) \\ 
		256 & 0.085 (0.98) & 6.77e-02 (0.98) & 7.58e-03 (1.00) \\ 
		512 & 0.185 (1.04) & 1.53e-01 (1.07) & 1.58e-02 (1.04) \\ 
		1024 & 0.384 (0.99) & 3.18e-01 (0.99) & 3.13e-02 (0.99) \\ 
		2048 & 0.816 (1.02) & 6.98e-01 (1.05) & 6.45e-02 (1.03) \\ 
		4096 & 1.736 (1.03) & 1.44e+00 (0.99) & 1.28e-01 (0.99)
	\end{tabular}
\end{table}

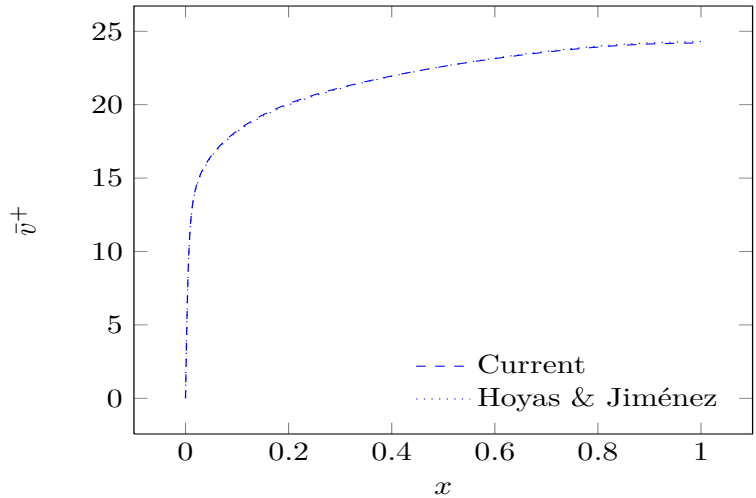
\begin{figure}[t]
	\begin{center}
		\begin{tikzpicture}[xscale=1.2]
		\begin{axis}[
		xlabel=$x$,
		ylabel=$\bar{v}^{+}$,
		legend cell align=left,
		legend pos=south east,
		legend style={fill=none, draw=none},
		]
		\addplot[dashed, color=blue] table [x=x, y=u]{my_data.txt};
		\addlegendentry{Current}
		\addplot[dotted, color=blue] table [x=x, y=u]{hoyas_data5.txt};
		\addlegendentry{Hoyas \& Jim{\' e}nez}
		\end{axis}
		\end{tikzpicture}
		\caption{Average velocity $\D{v}^+ = \D{v}/u_{\tau}$ in $y$-direction, for $Re_{\tau}=2000$. Our results (dashed) are compared with those of Hoyas and Jim\'{e}nez \cite{hoyas06} (dotted).}
		\label{fig:U_mean}
	\end{center}
\end{figure}

For the ultimate test of a fully turbulent flow, we use the turbulent Reynolds number $Re_{\tau}=2000$, which has been simulated previously by Hoyas and Jim\'{e}nez \cite{hoyas06}, Lozano and Jim\'{e}nez \cite{Lozano2014} and Lee and Moser \cite{leemoser15}. Lozano and Jim\'{e}nez showed that a small box  of size $2 \times 2\pi \times \pi$ in $x, y,$ and $z$ directions, respectively, could be simulated with $633 \times 1024 \times 1024$ collocation points (compact finite differences in wall-normal direction, Fourier in periodic) in order to capture all relevant one-point statistics found in much larger boxes. We use the same small computational box, with resolution of $512 \times 1024 \times 1024$, and expect the superior resolution properties of spectral methods to work in our favor. Furthermore, the driving force $\beta$ was dynamically adjusted in time to keep the flux through the channel constant. Simulations were run on the Shaheen II supercomputer at KAUST, using 512 cores. With a time step of $2\cdot10^{-5}$, the solver used approximately 1 second per time step. After reaching steady state, statistics were sampled for 100,000 time steps (enough to converge first order statistics), and the average velocity is shown in Figure \ref{fig:U_mean} in a dashed line, alongside the results obtained by Hoyas and Jim\'{e}nez \cite{hoyas06} in a dotted line. The difference is hardly visible, showing that the fully spectral Shen-Fourier solver is reliable and robust for $Re_{\tau} = 2000$. Furthermore, the results presented in Tables \ref{tab:roundoff}, \ref{tab:timings} and \ref{tab:timings_solver} indicate that the solver is more than likely to be reliable for even higher Reynolds numbers and larger scales.

\section{Conclusions}
In this paper a spectral-Galerkin Navier Stokes solver amenable to large-scale turbulent channel flows 
has been described. In the wall-normal direction, the solver utilizes base 
functions 
suggested by J. Shen \cite{Shen95}, constructed from Chebyshev polynomials. The 
solver uses Fourier decompositions in both periodic directions, and, as such, 
it is fully spectral in all spatial directions. To validate the method for large-scale simulations at high Reynolds numbers, we have shown 
that the roundoff errors are small, even for simulations twice the size of the 
largest simulations known to date \cite{leemoser15}.
The computational cost of the solver is shown to scale with the expected 
$\mathcal{O}(N_x \log N_x)$, for a mesh of size $N_x \times N_y \times N_z$ ($N_y$ and $N_z$ kept constant), due to fast 
transforms. Direct 
solvers have been devised for Helmholtz and biharmonic coefficient matrices, 
and these have been shown to scale close to 
$\mathcal{O}(N_x)$, with insignificant roundoff errors. Furthermore, all required scalar product matrices and fast transforms between spectral and physical space have been described. The fully spectral solver 
has been applied to a turbulent 
channel flow at $Re_{\tau}=2000$, and the first order statistics are shown to agree well with the literature. The solver has been implemented as part of the open source project spectralDNS \cite{spectralDNS}.

\section*{Acknowledgements}
This research is a part of the 4DSpace Strategic Research Initiative at the University of Oslo. I am also grateful to David Ketcheson and the KAUST Supercomputing Laboratory for providing access to Shaheen II.

\bibliography{bib.bib}

\section*{Appendix}
\label{sec:app}

\paragraph*{Fast scalar products and transforms}
Algorithms for fast scalar products and transforms, for the three 
spaces $W_N, \D{W}_N, \N{W}_N$, see Eqs. (\ref{eq:Tk}, \ref{eq:phiD}) and (\ref{eq:phiN}), are given in Alg. ~\ref{alg:fst}. The 
algorithms for the inverse transforms are given in Alg.~\ref{alg:ifst}. The 
algorithms are using the discrete cosine transforms (dct) of type 1, 2 or 3 
depending on the choice of discretization points and direction.

\paragraph*{Matrices}
Matrices representing one-dimensional scalar products are given in Table \ref{tab:matrices}. Note that the mass matrices are symmetrical, and as such only the upper triangular parts are described. Also note that the matrices are computed using quadrature, and not exact integration, which has some minor implications for Chebyshev-Gauss-Lobatto (GL) points. This follows since for an exact $L^2_{\sigma}$ scalar product we would have $c_{N_x}=1$ for GL as well as Chebyshev-Gauss (GC), but for the $l^2_{\sigma}$ scalar product used here $c_{N_x}=2$. This disagreement, that follows from inexact quadrature at the highest mode, explains the inclusion of the $c_{k+4}$ term for matrix $\N{B}$, and the $c_{k+2}$ term for matrix $\D{B}$. These are missing from Lemma 2.1 and 3.1 of \cite{Shen95}, because Shen here considers only the exact $L^2_{\sigma}$ scalar product.

\SetKwFunction{forwardScalar}{forwardScalar}
\SetKwFunction{forwardTransform}{forwardTransform}
\begin{algorithm}
    \DontPrintSemicolon
	\caption{Forward scalar products and transforms for all spaces $W_N, 
	\D{W}_N, \N{W}_N$. 
		Here "dct" is the discrete cosine transform from SciPy. The matrices $B$, $\D{B}$ and $\N{B}$ are given in Table \ref{tab:matrices}.\label{alg:fst}}
			\Fn{\forwardScalar{$\{f_j\}_{j=0}^{N}$, points, space}}{
            
            \Output{$\{y_k\}_{k=0}^{N}$}
                
			\If{ points = "Chebyshev-Gauss"}{
			$\{ w_k\}_{k=0}^{N} \gets \text{dct}(\{f_j\}_{j=0}^{N}, \text{type=2}, 
			\text{axis=0})\frac{\pi}{2 N}$ \;
            }
			\ElseIf{ points = "Chebyshev-Gauss-Lobatto"}{
			$\{ w_k\}_{k=0}^{N} \gets \text{dct}(\{f_j\}_{j=0}^{N}, \text{type=1}, 
			\text{axis=0}) \frac{\pi}{2 (N-1)}$ \;
		    }
			$\{y\}_{k=0}^N \gets 0$ \;
			\If{ space = $\D{W}_N$}{
			$\{y_k\}_{k=0}^{N-2} \gets \{w_k\}_{k=0}^{N-2} - \{w_k\}_{k=2}^{N}$ 
			\;
            }
			\ElseIf{ space = $\N{W}_N$}{
			\For{$k=0$ \KwTo $N-4$}{
			${y}_k \gets w_k - \frac{2(k+2)}{k+3} w_{k+2} + 
			\frac{k+1}{k+3}w_{k+4}$ \;
		    }
            }
			\ElseIf{space = $W_N$}{
			$\{y_k\}_{k=0}^N \gets \frac{2}{\pi} \{w_k\}_{k=0}^N $ \;
			$y_0 \gets y_0 / 2$ \;
			\If{ points = "Chebyshev-Gauss-Lobatto"}{
			$y_N \gets y_N / 2$ \;   	
		    }
            }
			\Return $\{y_k\}_{k=0}^{N}$    \;
            }
			\Fn{\forwardTransform{$\{f_j\}_{j=0}^{N}$, points, space}}{
			\Output{$\{y_k\}_{k=0}^{N}$}
            $\{y_k\}_{k=0}^{N} \gets$ forwardScalar($\{f_j\}_{j=0}^{N}$, 
\emph{points}, \emph{space}) \;
			\If{space = $\D{W}_N$}{
			$\{y_k\}_{k=0}^{N-2} \gets \{\D{B}_{kj}^{-1}\}_{k,j=0}^{N-2}\{y_k\}_{k=0}^{N-2}$ \;
            }
			\ElseIf{space = $\N{W}_N$}{
			$\{y_k\}_{k=0}^{N-4} \gets \{\N{B}_{kj}^{-1}\}_{k,j=0}^{N-4}\{y_k\}_{k=0}^{N-4}$ \;
            }
			\ElseIf{space = $W_N$}{
			$\{y_k\}_{k=0}^{N} \gets \{{B}_{kj}^{-1}\}_{k,j=0}^{N}\{y_k\}_{k=0}^{N}$ \;
            }
			\Return $\{y_k\}_{k=0}^{N}$
		}
\end{algorithm}
\SetKwFunction{inverseChebTransform}{inverseChebTransform}
\SetKwFunction{inverseShenTransform}{inverseShenTransform}
\begin{algorithm}
	\caption{Inverse transforms for all spaces $W_N, \D{W}_N, \N{W}_N$. 
		Here "dct" is the discrete cosine transform from SciPy.\label{alg:ifst}}
		\Fn {\inverseChebTransform{$\{f\}_{k=0}^N$, points}}{
        \Output{$\{y_j\}_{j=0}^{N}$}    
        \If{points = "Chebyshev-Gauss"}{
		    $\{y_j\}_{j=0}^{N} \gets \text{dct}(\{f\}_{k=0}^N, \text{type=3}, 
		    \text{axis=0})$ \;
		    $\{y_j\}_{j=0}^{N} \gets \frac{1}{2}\{y_j\}_{j=0}^{N} + 
		    \frac{1}{2}f_0$ \;
		}
        \ElseIf{ points = "Chebyshev-Gauss-Lobatto"}{
	        $\{y_j\}_{j=0}^{N} \gets \text{dct}(\{f\}_{k=0}^N, \text{type=1}, 
	        \text{axis=0})$ \;
	        $\{y_j\}_{j=0}^{N} \gets \frac{1}{2}\{y_k\}_{k=0}^{N} + 
	        \frac{1}{2}f_0$ \;
	        \For{$j=0$ \KwTo $N$}{
	         $y_j \gets y_j + \frac{(-1)^j}{2} f_N$ \;
    	    }
	    }
	    \Return{$\{y_j\}_{j=0}^{N}$}
    }
	    
	    \Fn{\inverseShenTransform{$\{f_k\}_{k=0}^{N}$, points, space}}{
		\Output{$\{y_j\}_{j=0}^{N}$}
        $\{z_k\}_{k=0}^{N} \gets 0$ \;
		\If{ space = $\D{W}_N$}{
		$\{z_k\}_{k=0}^{N-2} \gets \{f_k\}_{k=0}^{N-2}$ \;
		$\{z_k\}_{k=2}^{N} \gets \{z_k\}_{k=2}^{N} - \{f_k\}_{k=0}^{N-2} $ \;
        }
        \ElseIf{ space = $\N{W}_N$}{	       
		$\{z_k\}_{k=0}^{N-4} \gets \{f_k\}_{k=0}^{N-4}$ \;
		\For{$k=0$ \KwTo $N-4$}{
		$z_{k+2} \gets z_{k+2} - \frac{2(k+2)}{k+3} f_k$ \;
	    }
		\For{$k=0$ \KwTo $N-4$}{
		$z_{k+4} \gets z_{k+4} + \frac{k+1}{k+3} f_k$ \;
	    }
        }
        $\{y_j\}_{j=0}^{N} \gets $ inverseChebTransform($\{z_k\}_{k=0}^N$, 
        \emph{points}) 
        \;
		\Return $\{y_j\}_{j=0}^{N}$ \;
    }
\end{algorithm}
\begin{table}
	\centering
	\caption{Matrix notation and description for one dimensional scalar products. Note that $c_0=2$ and $c_k=1$ for $0<k<N_x$. For Chebyshev-Gauss points $c_{N_x}=1$, whereas $c_{N_x}=2$ for Chebyshev-Gauss-Lobatto points. Also note that only the upper diagonals are described for the symmetric matrices.	\label{tab:matrices}}
	\begin{tabular}{ccl}	
		Notation & Scalar product & Description \\ 
		\hline
%

$\{\N{C}_{kj}\}_{k,j=0}^{N_x-4, N_x-2}$ & $\left(\D{\phi}^{'}_j, \N{\phi}_k 
\right)_{\sigma}$ & $\begin{cases} -\pi(k+1), &j=k-1,\\
2\pi(k+1), & j=k+1, \\
-\pi(k+1), & j=k+3, \\
0, &\text{otherwise.} \end{cases}$ \\

$\{\D{C}_{kj}\}_{k,j=0}^{N_x-2, N_x-4}$ & $\left(\N{\phi}^{'}_j, \D{\phi}_k 
\right)_{\sigma}$ & $\begin{cases}
\pi (k-2)(k+1)/k, &j=k-3,\\
-2 \pi \frac{(k+1)^2}{k+2}, & j=k-1, \\
\pi(k+1), & j=k+1, \\
0, &\text{otherwise.}
\end{cases}$ \\

$\{{C}_{kj}\}_{k,j=0}^{N_x, N_x-2}$ & $\left(\D{\phi}^{'}_j, T_k 
\right)_{\sigma}$ & $\begin{cases}
-\pi (k+1), &j=k-1,\\
-2 \pi, & j=k+1,k+3, \ldots, N_x-2 \\
0, &\text{otherwise.}
\end{cases}$ \\

$\{\N{M}_{kj}\}_{k,j=0}^{N_x-4, N_x-2}$  & $ \left(\D{\phi}_j, \N{\phi}_k \right)_{\sigma} $ & $\begin{cases}
-\pi/2 & j=k-2, \\
\frac{\pi}{2} \left(c_k+2 \frac{k+2}{k+3} \right) & j=k, \\
-\frac{\pi}{2} \left(2\frac{k+2}{k+3} + c_{j+2} \frac{k+1}{k+3}\right) & j=k+2 \\
\frac{\pi}{2}\frac{k+1}{k+3} & j=k+4, \\
0, &\text{otherwise.}
\end{cases}$ \\

$\{\N{B}_{kj}\}_{k,j=0}^{N_x-4}$ & $(\N{\phi}_j, \N{\phi}_k)_{\sigma}$ & $\begin{cases}
\frac{\pi}{2} \left(c_k + 4 \left(\frac{k+2}{k+3} \right)^2 + c_{k+4} 
\left(\frac{k+1}{k+3}\right)^2    \right), &j=k,\\
-\pi \left( \frac{k+2}{k+3} + \frac{(k+4)(k+1)}{(k+5)(k+3)} \right), &j=k + 2,\\
\frac{\pi}{2} \frac{k+1}{k+3} , & j=k + 4, \\
0, &\text{otherwise.}
\end{cases} $\\

$\{\D{B}_{kj}\}_{k,j=0}^{N_x-2} $ & $(\D{\phi}_j, \D{\phi}_k)_{\sigma}$ & $ \begin{cases} 
\frac{\pi}{2} (c_k+c_{k+2}) &j=k, \\
-\frac{\pi}{2} &j=k + 2, \\
0, &\text{otherwise.}
\end{cases}$ \\

$\{B_{kj}\}_{k,j=0}^{N_x} $ & $(T_k, T_j)_{\sigma}$ & $\frac{c_k \pi}{2} \delta_{kj}.$

	\end{tabular}
\end{table}

\end{document}